\documentclass[oneside,english]{amsart}
\usepackage[T1]{fontenc}
\usepackage[latin9]{inputenc}
\usepackage{amsthm}
\usepackage{amstext}
\usepackage{amssymb}

\makeatletter
\numberwithin{equation}{section}
\numberwithin{figure}{section}
\theoremstyle{plain}
\newtheorem{thm}{\protect\theoremname}
  \theoremstyle{plain}
  \newtheorem{prop}[thm]{\protect\propositionname}

\usepackage{ae,aecompl}

\usepackage{fullpage}

\usepackage{pifont}

\makeatother

\usepackage{babel}
  \providecommand{\propositionname}{Proposition}
\providecommand{\theoremname}{Theorem}

\begin{document}
\global\long\def\CC{\mathbb{C}}

\global\long\def\RR{\mathbb{R}}

\global\long\def\EE{\mathcal{E}}

\global\long\def\FF{\mathcal{F}}

\global\long\def\phi{\varphi}

\global\long\def\half{\nicefrac{1}{2}}

\global\long\def\re{\operatorname{Re}}

\title{A letter: The log-Brunn-Minkowski inequality for complex bodies}

\author{Liran Rotem}

\maketitle
We will use the following terminology: A real body $K\subseteq\RR^{n}$
is the unit ball of a norm $\left\Vert \cdot\right\Vert $ on $\RR^{n}$,
i.e. a convex, origin symmetric, compact set with non-empty interior.
Similarly, a complex body $K\subseteq\CC^{n}$ is the unit ball of
a norm $\left\Vert \cdot\right\Vert $ on $\CC^{n}$. By identifying
$\CC^{n}\simeq\RR^{2n}$ we see that every complex body is also a
real body, but not vice versa. In fact, a complex body $K\subseteq\CC^{n}$
is a real body which is also symmetric with respect to complex rotations,
i.e. if $z\in K$ implies that $e^{i\theta}z\in K$ for all $\theta\in\RR$. 

For a real body $K$, the support function of $K$ is defined as $h_{K}(\theta)=\left\Vert \theta\right\Vert _{K}^{\ast}=\sup_{x\in K}\left\langle x,\theta\right\rangle $.
Given two such bodies $K$ and $T$ and a number $0\le\lambda\le1$,
we define the logarithmic mean of $K$ and $T$ by 
\[
L_{\lambda}(K,T)=\left\{ x\in\RR^{n}:\ \left\langle x,\theta\right\rangle \le h_{K}(\theta)^{1-\lambda}h_{T}(\theta)^{\lambda}\text{ for all }\theta\in\RR^{n}\right\} ,
\]
where $\left\langle \cdot,\cdot\right\rangle $ denotes the standard
Euclidean inner product. 

The log-Brunn-Minkowski inequality states that $\left|L_{\lambda}(K,T)\right|\ge\left|K\right|^{1-\lambda}\left|T\right|^{\lambda}$,
where $\left|\cdot\right|$ denotes the (Lebesgue) volume. It was
conjectured by Böröczky, Lutwak, Yang and Zhang (\cite{Boroczky2012}),
who proved it for $K,T\subseteq\RR^{2}$. Saroglou proved (\cite{Saroglou2014})
that the inequality holds when $K$ and $T$ are $n$-dimensional
real bodies which are unconditional with respect to the same basis.

The goal of this note is to explain why the log-Brunn-Minkowski inequality
holds for complex bodies:
\begin{thm}
\label{thm:main-thm}For complex bodies $K,T\subseteq\CC^{n}$ and
$0\le\lambda\le1$ we have $\left|L_{\lambda}(K,T)\right|\ge\left|K\right|^{1-\lambda}\left|T\right|^{\lambda}$.
\end{thm}
Theorem \ref{thm:main-thm} will follow from a result of Cordero-Erausquin
(\cite{Cordero-Erausquin2002}). In his work, Cordero-Erausquin proved
a generalization of the Blaschke-Santaló inequality in the complex
case. Specifically, he proved that for complex bodies $K,T\subseteq\CC^{n}$
we have 
\begin{equation}
\left|K\cap T\right|\left|K^{\circ}\cap T\right|\le\left|B_{2}^{2n}\cap T\right|,\tag{\ensuremath{\ast}}\label{eq:gen-santalo}
\end{equation}
 where $K^{\circ}$ is the polar body to $K$ and $B_{2}^{2n}\subseteq\CC^{n}$
is the unit Euclidean ball. As a side note we remark that proving
the same inequality for general real bodies is an open problem --
see \cite{Klartag2007c} for a partial result and a short discussion.

Cordero-Erausquin's proved the inequality (\ref{eq:gen-santalo})
as a corollary of a general theorem about complex interpolation -
see Theorem \ref{thm:complex-lc} below. The main point of this letter
is the observation that the same general theorem also implies Theorem
\ref{thm:main-thm}. This was apparently known to Cordero-Erausquin
himself, but not to other researchers in the community who haven't
studied the complex case carefully. Theorem \ref{thm:main-thm} may
be a strong indication that the log-Brunn-Minkowski conjecture is
true in general. Alternatively, it may indicate the existence of a
rich theory of geometric inequalities in the complex case.

Let us briefly recall the definition of complex interpolation. We
will give the construction for the finite-dimensional case, following
the presentation of \cite{Cordero-Erausquin2002}, and refer the reader
to \cite{Bergh1976} for a more detailed account. Set $S=\left\{ z\in\CC:\ 0<\re z<1\right\} $,
and define 
\[
\FF=\left\{ f:\overline{S}\to\CC^{n}:\ \begin{array}{l}
f\text{ is bounded and continuous on \ensuremath{\bar{S}} and analytic on }S\\
\text{such that }{\displaystyle \lim_{t\to\pm\infty}f(it)=\lim_{t\to\pm\infty}f(1+it)=0}
\end{array}\right\} .
\]
 Given two norms $\left\Vert \cdot\right\Vert _{0}$ and $\left\Vert \cdot\right\Vert _{1}$
on $\CC^{n}$, we define a norm on $\FF$ by 
\[
\left\Vert f\right\Vert _{\FF}=\max\left\{ \sup_{t\in\RR}\left\Vert f(it)\right\Vert _{0},\sup_{t\in\RR}\left\Vert f(1+it)\right\Vert _{1}\right\} .
\]
Finally, for $\lambda\in[0,1]$, we define the interpolated norm $\left\Vert \cdot\right\Vert _{\lambda}$
by 
\[
\left\Vert x\right\Vert _{\lambda}=\inf\left\{ \left\Vert f\right\Vert _{\mathcal{F}}:\ f\in\FF,\ f(\lambda)=x\right\} .
\]

It is not hard to see that for $\lambda=0,1$ we recover the original
norms $\left\Vert \cdot\right\Vert _{0},\left\Vert \cdot\right\Vert _{1}$.
The only other result we will need from the standard theory of complex
interpolation is the following:
\begin{prop}
\label{prop:interp-bound}Let $\left\Vert \cdot\right\Vert _{0}$,\textup{
$\left\Vert \cdot\right\Vert _{1}$ be norms on $\CC^{n}$ and let
$\left\Vert \cdot\right\Vert _{\lambda}$ be the interpolated norms.
Then 
\[
\left\Vert z\right\Vert _{\lambda}^{\ast}\le\left(\left\Vert z\right\Vert _{0}^{\ast}\right)^{1-\lambda}\left(\left\Vert z\right\Vert _{1}^{\ast}\right)^{\lambda}
\]
 for every $z\in\CC^{n}.$}
\end{prop}
This inequality, with its simple proof, may be found for example in
\cite{Pisier1989} as equation $\left(7.26\right)^{\ast}$. 

If $K$ is the unit ball of $\left\Vert \cdot\right\Vert _{0}$ and
$T$ is the unit ball of $\left\Vert \cdot\right\Vert _{1}$ we will
write $C_{\lambda}(K,T)$ for the unit ball of $\left\Vert \cdot\right\Vert _{\lambda}$.
Proposition \ref{prop:interp-bound} implies that $h_{C_{\lambda}(K,T)}(z)\le h_{K}(z)^{1-\lambda}h_{T}(z)^{\lambda}$
for all $z\in\CC^{n}$, and hence $C_{\lambda}(K,T)\subseteq L_{\lambda}(K,T)$.

The main theorem of \cite{Cordero-Erausquin2002} is the following:
\begin{thm}
\label{thm:complex-lc}The function $\lambda\longmapsto\left|C_{\lambda}(K,T)\right|$
is log-concave on $[0,1]$. 
\end{thm}
It is now easy to deduce Theorem \ref{thm:main-thm}, as 
\[
\left|L_{\lambda}(K,T)\right|\ge\left|C_{\lambda}(K,T)\right|\ge\left|C_{0}(K,T)\right|^{1-\lambda}\cdot\left|C_{1}(K,T)\right|^{\lambda}=\left|K\right|^{1-\lambda}\left|T\right|^{\lambda}.
\]

\bibliographystyle{plain}
\bibliography{C:/Users/Liran/Dropbox/citations/library}

\end{document}